# Hasse-Schmidt derivations
# and the
# Hopf algebra of non-commutative symmetric functions

Michiel Hazewinkel
Burg. 's Jacob laan  18
NL-1401BR  BUSSUM
The Netherlands
<michhaz@xs4all.nl>

1. **Introduction**.
Let $A$ be an associative algebra (or any other kind of algebra for that matter). A derivation on $A$ is an endomorphism $\partial$ of the underlying Abelian group of $A$ such that

$$\partial(ab) = a(\partial b) + (\partial a)b \quad \text{for all} \quad a,b \in A \tag{1.1}$$

A Hasse-Schmidt derivation is a sequence $(d_0 = \mathrm{id}, d_1, d_2, \cdots, d_n, \cdots)$ of endomorphisms of the underlying Abelian group such that for all $n \geq 1$

$$d_n(ab) = \sum_{i=0}^{n} (d_i a)(d_{n-i} b) \tag{1.2}$$

Note that $d_1$ is a derivation as defined by (1.1). The individual $d_n$ that occur in a Hasse-Schmidt derivation are also sometimes called higher derivations.

A question of some importance is whether Hasse-Schmidt derivations can be written down in terms of polynomials in ordinary derivations. For instance in connection with automatic continuity for Hasse-Schmidt derivations on Banach algebras.

Such formulas have been written down by, for instance, Heerema and Mirzavaziri in [5] and [6]. They also will be explicitly given below.

It is the purpose of this short note to show that such formulas follow directly from some easy results about the Hopf algebra **NSymm** of non-commutative symmetric functions. In fact this Hopf algebra constitutes a universal example concerning the matter.

2. **Hopf algebras and Hopf module algebras**. Everything will take place over a commutative associative unital base ring $k$; unadorned tensor products will be tensor products over $k$. In this note $k$ will be the ring of integers **Z**, or the field of rational numbers **Q**.

Recall that a Hopf algebra over $k$ is a $k$-module $H$ together with five $k$-module morphisms

$m: H \otimes H \longrightarrow H$, $e: k \longrightarrow H$, $\mu: H \longrightarrow H \otimes H$, $\varepsilon: H \longrightarrow k$, $\iota: H \longrightarrow H$ such that $(H, m, e)$ is an associative $k$-algebra with unit, $(H, \mu, \varepsilon)$ is a co-associative co-algebra with co-unit, $\mu$ and $\varepsilon$ are algebra morphisms (or, equivalently, that $m$ and $e$ are co-algebra morphisms), and such that $\iota$ satisfies $m(\iota \otimes \mathrm{id})\mu = \varepsilon e$, $m(\mathrm{id} \otimes \iota)\mu = \varepsilon e$. The antipode $\iota$ will play no role in what follows. If there is no antipode (specified) one speaks of a bi-algebra. For a brief introduction to Hopf algebras (and co-algebras) with plenty examples see e.g. chapters 2 and 3 of [4].

Recall also that an element $p \in H$ is called primitive if $\mu(p) = p \otimes 1 + 1 \otimes p$. These form a sub-$k$-module of $H$ and form a Lie algebra under the commutator difference product $(p, p') \mapsto pp' - p'p$. I shall use $Prim(H)$ to denote this $k$-Lie-algebra.

Given a Hopf algebra over $k$, a Hopf module algebra is a $k$-algebra $A$ together with an action of the underlying algebra of $H$ on (the underlying module of) $A$ such that moreover

$$h(ab) = \sum_{(h)} h_{(1)}(a) h_{(2)}(b) \text{ for all } a, b \in A, \text{ where } \mu(h) = \sum_{(h)} h_{(1)} \otimes h_{(2)} \tag{2.1}$$

where I have used to Sweedler-Heynemann notation for the co-product.

Note that this means that the primitive elements of $H$ act as derivations.

3. **The Hopf algebra NSymm of non-commutative symmetric functions**. As an algebra over the integers **NSymm** is simply the free associative algebra in countably many (non-commuting) indeterminates, $\mathbf{NSymm} = \mathbf{Z}\langle Z \rangle = \mathbf{Z}\langle Z_1, Z_2, \cdots \rangle$. The comultiplication and counit are given by

$$\mu(Z_n) = \sum_{i+j=n} Z_i \otimes Z_j, \text{ where } Z_0 = 1, \; \varepsilon(1) = 1, \; \varepsilon(Z_n) = 0 \text{ for } n \geq 1 \tag{3.1}$$

As **NSymm** is free as an associative algebra it is no trouble to verify that his defines a bi-algebra. The seminal paper [1] started the whole business of non-commutative symmetric functions, now a full-fledged research area in its own right.

Now consider an **NSymm** Hopf module algebra $A$. Then, by (2.1) and (3.1) the module endomorphims defined by by the actions of the $Z_n$, $n \geq 1$, $d_n(a) = Z_n a$, define a Hasse-Schmidt derivation. Inversely, if $A$ is a $k$-algebra together with a Hasse-Schmidt derivation one defines a **NSymm** Hopf module algebra structure on $A$ by setting $Z_n a = d_n(a)$. This works because **NSymm** is free as an algebra.

Thus an **NSymm** Hopf module algebra $A$ is precisely the same thing as a $k$-algebra $A$ together with a Hasse-Schmidt derivation on it and the matter of writing the elements of the sequence of morphisms that make up the Hasse-Schmidt derivation in terms of ordinary derivations comes down to the matter of finding enough primitives of **NSymm** so that the generators, $Z_n$, can be written as polynomials in these primitives.

4. **The Newton primitives of NSymm**. Define the non-commutative polynomials $P_n$ and $P_n'$ by the recursion formulas

$$\begin{aligned} P_n &= nZ_n - (Z_{n-1}P_1 + Z_{n-2}P_2 + \cdots + Z_1 P_{n-1}) \\ P_n' &= nZ_n - (P_1' Z_{n-1} + P_2' Z_{n-2} + \cdots + P_{n-1}' Z_1) \end{aligned} \tag{4.1}$$

These are non-commutative analogues of the well known Newton formulas for the power sums in terms of the complete symmetric functions in the usual commutative theory of symmetric functions. It is not difficult to write down an explicit expression for these polynomials:

$$P_n(Z) = \sum_{\substack{i_1 + \cdots + i_m = n \\ i_j \in \mathbf{N}}} (-1)^{m+1} i_m Z_{i_1} Z_{i_2} \cdots Z_{i_m} \tag{4.2}$$

Nor is it difficult to write down a formula for the $Z_n$ in terms of the $P$'s or $P''$s.[1] However, to do that one definitely needs to use rational numbers and not just integers. For instance

$$Z_2 = \frac{P_1^2 + P_2}{2}$$

The key observation is now:

4.3. **Proposition**. The elements $P_n$ and $P_n'$ are primitive elements of the Hopf algebra **NSymm**.

The proof is a straightforward uncomplicated induction argument using the recursion formulas (4.1). See e.g. [4], p. 147.

Using the $P_n'$ an immediate corollary is the following main theorem from [6].

4.4. **Theorem**. Let $A$ be an associative algebra over the rational numbers $\mathbf{Q}$ and let $(\mathrm{id}, d_1, d_2, \cdots, d_n, \cdots)$ be a Hasse-Schmidt derivation on it. Then the $\delta_n$ defined recursively by

$$\delta_n = n d_n - \delta_1 d_{n-1} - \cdots - \delta_{n-1} d_1 \tag{4.5}$$

are ordinary derivations and

$$d_n = \sum_{\substack{r_1 + r_2 + \cdots r_m = n \\ r_j \in \mathbf{N}}} c_{r_1, r_2, \cdots, r_m} \delta_{r_1} \delta_{r_2} \cdots \delta_{r_m}, \tag{4.6}$$

where

$$c_{r_1, r_2, \cdots, r_m} = \frac{1}{r_1 + r_2 + \cdots + r_m} \frac{1}{r_2 + \cdots + r_m} \cdots \frac{1}{r_{m-1} + r_m} \frac{1}{r_m} \tag{4.7}$$

4.8. **Comment**. Because

$$P_n' \equiv n Z_n \mod(Z_1, Z_2, \cdots, Z_{n-1})$$

---

[1] This is an instance where the noncommutative formulas are are more elegant and also easier to prove than their commutative analogues. In the commutative case there are all kinds of multiplicities that mess things up.

the formulas expressing the $Z_n$ in terms of the $P'_n$ are unique and so denominators are really needed.

4.9. **Comment and example**. There are many more primitive elements in **NSymm** than just the $P'_n$ and $P_n$. One could hope that by using all of them integral formulas for the $Z_n$ in terms of primitives would become possible. That is not the case. The full Lie algebra of primitives of **NSymm** was calculated in [3]. It readily follows from the description there that $\mathbf{Z}\langle Prim(\mathbf{NSymm})\rangle$, the sub-algebra of **NSymm** generated by all primitive elements is strictly smaller than **NSymm**. In fact much smaller in a sense that can be specified. Thus the theorem does not hold over the integers.

A concrete example of a Hasse-Schmidt derivation of which the constituting endomorphisms can not be written as integral polynomials in derivations can be given in terms of **NSymm** itself. As follows. The Hopf algebra **NSymm** is graded by giving $Z_n$ degree $n$. Note that each graded piece is a free **Z**-module of finite rank. Let **QSymm**, often called the Hopf algebra of quasi-symmetric functions, be the graded dual Hopf algebra. Then each $Z_n$ defines a functional $\alpha_n : \mathbf{QSymm} \longrightarrow \mathbf{Z}$. Now define an endomorphism $d_n$ of **QSymm** as the composed morphism

$$\mathbf{QSymm} \xrightarrow{\mu_{\mathbf{QSymm}}} \mathbf{QSymm} \otimes_{\mathbf{Z}} \mathbf{QSymm} \xrightarrow{id \otimes \alpha_n} \mathbf{QSymm}$$

Then the $d_n$ form a Hasse-Schmidt derivation of which the components can not be written as integer polynomials in ordinary derivations.

5. **The Hopf algebra LieHopf**.

In [5] there occurs a formula for manufacturing Hasse-Schmidt derivations from a collection of ordinary derivations that is more pleasing – at least to me – than 4.6. This result from loc. cit. can be strengthened to give a theorem similar to theorem 4.4 but with more symmetric formulae. This involves another Hopf algebra over the integers which I like to call **LieHopf**.

As an algebra **LieHopf** is again the free associative algebra in countably many indeterminates $\mathbf{Z}\langle U \rangle = \mathbf{Z}\langle U_1, U_2, \cdots \rangle$. But this time the co-multiplication and co-unit are defined by

$$\mu(U_n) = U_n \otimes 1 + 1 \otimes U_n, \quad \varepsilon(U_n) = 0 \tag{5.1}$$

so that all the $U_n$ are primitive. And, in fact the Lie algebra of primitives of this Hopf algebra is the free Lie algebra on countably many generators.

Over the integers **LieHopf** and **NSymm** are very different but over the rationals they become isomorphic. There are very many isomorphisms. A particularly nice one is given by considering the power series identity

$$1 + Z_1 t + Z_2 t^2 + Z_3 t^3 + \cdots = \exp(U_1 t + U_2 t^2 + U_3 t^3 + \cdots) \tag{5.2}$$

which gives the following formulae for the $U$'s in terms of the $Z$'s and vice versa.

$$Z_n(U) = \sum_{r_1 + \cdots + r_m = n} \frac{U_{r_1} U_{r_2} \cdots U_{r_m}}{m!} \tag{5.3}$$

$$U_n(Z) = \sum_{r_1+\cdots+r_m=n} (-1)^{m+1} \frac{Z_{r_1}Z_{r_2}\cdots Z_{r_m}}{m} \tag{5.4}$$

For two detailed proofs that these formulas do indeed give an isomorphism of Hopf algebras see [2]; or see chapter 6 of [4]. In terms of derivations, reasoning as above, this gives the following theorem.

5.5. **Theorem**. Let $A$ be an algebra over the rationals and let $(\mathrm{id},d_1,d_2,\cdots)$ be a Hasse-Schmidt derivation on it. Then the $\partial_n$ defined by

$$\partial_n = \sum_{r_1+\cdots+r_m=n} (-1)^{m+1} \frac{d_{r_1}d_{r_2}\cdots d_{r_m}}{m} \tag{5.6}$$

are (ordinary) derivations and

$$d_n = \sum_{r_1+\cdots+r_m=n} \frac{\partial_{r_1}\partial_{r_2}\cdots \partial_{r_m}}{m!} \tag{5.7}$$

5.8. **Comment**. Perhaps I should add that given any collection of ordinary derivations, formula (5.7) yields a Hasse-Schmidt derivation. That is the theorem from [5] with which I started this section.